\documentclass[twoside,11pt,reqno]{amsart}

\usepackage{upref,amsxtra,amssymb,amscd}
\usepackage{varioref}
\usepackage{verbatim}
\usepackage{epsfig}
\usepackage{color}

\usepackage{epic,eepic,eucal}
\usepackage{enumerate}

\def\Ba{        \textrm{\bf Bad}}

\def\a{         \alpha}
\def\b{         \beta}

\newcommand{\NN}{{\mathbb N}}
\newcommand{\RR}{{\mathbb R}}

\newcommand{\ZZ}{{\mathbb Z}}

\newtheorem{theo}{\sc Theorem}[section]

\newtheorem{lemm}[theo]{\sc Lemma}
\newtheorem{coro}[theo]{\sc Corollary}

\theoremstyle{definition}

\theoremstyle{remark}

\newtheorem{rema}[theo]{\sc Remark}

\numberwithin{equation}{section}

\begin{document}
\title{Badly approximable systems of affine forms, fractals, and Schmidt games}
\author{Manfred Einsiedler}
\author{Jimmy Tseng}
\address{Manfred Einsiedler, Department of Mathematics, The Ohio State University, Columbus, OH 43210}
\email{manfred@math.ohio-state.edu}
\address{Jimmy Tseng, Department of Mathematics, The Ohio State University, Columbus, OH 43210}
\email{tseng@math.ohio-state.edu}
\thanks{M.E. acknowledges the support of the NSF (DMS-grant 0554373) and of the SNF
(200021-127145).}

\begin{abstract}
A badly approximable system of affine forms is determined by a matrix and a vector.  We show Kleinbock's conjecture for badly approximable systems of affine forms:  for any fixed vector, the set of badly approximable systems of affine forms is winning (in the sense of Schmidt games) even when restricted to a fractal (from a certain large class of fractals).  In addition, we consider fixing the matrix instead of the vector where an analog statement holds.  

\end{abstract}

\maketitle
\section{Introduction}\label{secIntro}
Let $M_{m,n}(\RR)$ denote  the set of $m \times n$ real matrices and let $\widetilde{M}_{m,n}(\RR)$ denote $M_{m,n}(\RR) \times \RR^m$.  The element in $\widetilde{M}_{m,n}(\RR)$ corresponding to $A \in M_{m,n}(\RR)$ and $\mathbf{b} \in \RR^m$ will be expressed as $\langle A, \mathbf{b} \rangle$.  Consider the following well-known sets from the theory of Diophantine approximation (or metric number theory), see for instance \cite{Kl}:  
\begin{multline*}
 \Ba(m,n) := \Big\{\langle A, \mathbf{b} \rangle \in \widetilde{M}_{m,n}(\RR)  \mid \mbox{ there exists }  c(A, \mathbf{b}) > 0\\ \mbox{ such that } \|A \mathbf{q} - \mathbf{b} \|_{\ZZ} \geq \frac{c(A, \mathbf{b})}{\|\mathbf{q}\|^{n/m}} \mbox{ for all } \mathbf{q} \in \ZZ^n \backslash \{\mathbf{0}\}\Big\}
\end{multline*}
where $\|\cdot\|$ is the sup norm on $\RR^k$ and $\| \cdot \|_{\ZZ}$ is the function on $\RR^k$ given by $\|\mathbf{x}\|_{\ZZ} := \inf_{p \in \ZZ^k}\|\mathbf{x} - \mathbf{p}\|$.  The set $\Ba(m,n)$ is called the \textit{set of badly approximable systems of $m$ affine forms in $n$ variables.}  For any $\mathbf{b} \in \RR^m$, let $\Ba^{\mathbf{b}}(m,n) := \{A \in M_{m,n}(\RR) \mid \langle A, \mathbf{b} \rangle \in \Ba(m,n)\}$, and, for any $A \in M_{m,n}(\RR)$, let $\Ba_A(m,n):= \{\mathbf{b} \in \RR^m \mid  \langle A, \mathbf{b} \rangle \in \Ba(m,n)\}$.

The set $\Ba^{\mathbf{0}}(m,n)$ is called the \textit{set of badly approximable systems of $m$ linear forms in $n$ variables} and is an important and classical object of study in metric number theory.  Although $\Ba^{\mathbf{0}}(m,n)$ is a Lebesgue null set (Khintchine, 1926), it has full Hausdorff dimension and, even stronger, is winning as shown by Schmidt \cite{Sch3} in 1969.\footnote{One can even intersect $\Ba^{\mathbf{0}}(m,n)$ with certain fractals and still retain the winning property, see Theorem~1 of~\cite{F1}.}  Winning sets have a few other properties besides having full Hausdorff dimension.  An important example of such is the countable intersection property, which allows countable intersections of winning sets to remain winning.  This puts the class of winning sets next to other important classes of large sets with the same property as, for example, the class of conull sets or the class of dense $G_\delta$-sets.  In contrast, the class of sets that are simply of full Hausdorff dimension does not have the countable (or even finite) intersection property.  See Section~\ref{subsecBWSCF} for more details on the properties of winning sets. 

For general $\mathbf{b}$, less has heretofore been known.  Another result of Schmidt implies that $\Ba^{\mathbf{b}}(m,n)$ has zero Lebesgue measure for any $\mathbf{b}$~\cite{Sch1}.  With regard to dimension, however, D.~Kleinbock has shown that $\Ba^{\mathbf{b}}(m,n)$ has full Hausdorff dimension for $\mathbf{b}$ from a full Hausdorff dimension subset of $\RR^m$~\cite{Kl}.  Thus, a fundamental question in the theory of badly approximable systems of affine forms (and in metric number theory) is whether $\Ba^{\mathbf{b}}(m,n)$ has full Hausdorff dimension for every $\mathbf{b}$.  In fact, Kleinbock \cite{Kl} conjectured  that $\Ba^{\mathbf{b}}(m,n)$ is winning for every $\mathbf{b}$.  In this paper, we show that Kleinbock's conjecture is true and, moreover, that $\Ba^{\mathbf{b}}(m,n)$ is  winning even when restricted to certain fractals; see Theorem~\ref{thmBAwinning}.

Recently, interest in the size of related sets, namely the size of $\Ba_A(m,n)$ for fixed $A$, has developed.\footnote{Problems in metric number theory in which the vector $\mathbf{b}$ is fixed are referred to as \textit{singly metric inhomogeneous problems}.  Problems in which nothing is fixed are referred to as \textit{doubly metric inhomogeneous problems}.  Problems in which the matrix $A$ is fixed first appeared in this generality, the authors believe, in~\cite{BHKV} and are not, as of yet, named.}  The sets $\Ba_A(m,n)$ naturally arise as the complements of sets of toral translation vectors that satisfy certain shrinking target properties (see~\cite{T1} and~\cite{BHKV} for details).  For almost every $A$, these sets are Lebesgue null sets, but it is easy to see that these sets can possibly have even full Lebesgue measure.  However, regardless of Lebesgue measure, Y.~Bugeaud, S.~Harrap, S.~Kristensen, and S.~Velani have recently shown that, for every $A$, $\Ba_A(m,n)$ has full Hausdorff dimension even when restricted to certain fractals (Theorem~2 of~\cite{BHKV}).  Two questions are inspired by their result:  are the sets $\Ba_A(m,n)$ winning for all $A$ and, if so, can this winning property be further generalized to fractals from a larger class of fractals than those considered in~\cite{BHKV}.  In~\cite{T2}, the second-named author has answered the first question in the affirmative for the special case of $n=m=1$.  In this paper, we answer both questions in the affirmative for the general case; see Theorem~\ref{thmDualBAwinning}.

Finally, as a corollary of the proof of Theorem~\ref{thmDualBAwinning}, we also study the set of \textit{infinitely badly approximable matrices} 
\[
\Ba^\infty_A(m,n) := \Big\{\mathbf{b} \in \RR^m  \mid \liminf_{\mathbf{q} \in \ZZ^n \backslash \{\mathbf{0}\}} \|\mathbf{q}\|^{n}\|A \mathbf{q} - \mathbf{b} \|^m_{\ZZ} = \infty \Big\}
\] 
for matrices $A$ that are singular (in the sense of the theory of Diophantine approximation). Here we say that $A$ is \textit{singular} if for every $\varepsilon >0$ and large enough $N$ there are solutions  $\mathbf{q} \in \ZZ^n$ to the system of inequalities
\[
  \|A \mathbf{q}\|_{\ZZ} \leq \frac{\varepsilon}{N^{n/m}}\textrm{ and } 0 < \|\mathbf{q}\| < N.
\] 
We note that $\Ba^\infty_A(m,n) \subset \Ba_A(m,n)$.
 The set of singular matrices $A$, which we denote by $\mathbf{SM}_{m,n}(\RR)$, is called \textit{the set of singular systems of $m$ linear forms in $n$ variables} (or \textit{the set of singular $\mathbf{m \times n}$ matrices}) and is another important and classical object of study in metric number theory.

\subsection{Statement of results}\label{subsecStateResults} 

In this section, we state and discuss our results.  Note that $\dim(\cdot)$ refers to Hausdorff dimension throughout this paper and $\underline{d}_\mu(U)$ refers to lower pointwise dimension.\footnote{Recall that, for an open set $U$ of a metric space with a locally finite Borel measure $\mu$, the lower pointwise dimension is defined as \[\underline{d}_\mu (U):= \inf_{x \in U} \liminf_{r \rightarrow 0} \frac {\log \mu(B(x,r))}{\log r}.\]}  Our first result, Theorem~\ref{thmBAwinning}, answers affirmatively the aforementioned fundamental question in the theory of badly approximable systems of affine forms and, moreover, subsumes both the classical theory concerning the size of $\Ba^{\mathbf{0}}(m,n)$, which culminated in Schmidt's proof of the winning property, and the more recent proofs of L.~Fishman (\cite{F1} and~\cite{F2}) involving  the intersection of $\Ba^{\mathbf{0}}(m,n)$ with certain fractals. 

\begin{theo}\label{thmBAwinning}
Let $K \subset M_{m,n}(\RR)$ be the support of an absolutely friendly measure $\mu$ (as an example, the Lebesgue measure restricted to $[0,1]^{mn}$).  Then, for any $\mathbf{b} \in  \RR^m$, \[K \cap \Ba^{\mathbf{b}}(m,n)\] is a winning set on $K$.
\end{theo}

Note that the notions of winning and absolute decaying are defined in Section~\ref{secBackground}, but we note that the Lebesgue measure on $M_{m,n}(\RR)$ is absolutely decaying.  Also, we show that the winning parameter (see Section~\ref{secBackground} for the definition) is independent of $\mathbf{b}$.  For its value, see the proof of the result in Section~\ref{secAPothmBAwinning}.

Theorem~\ref{thmBAwinning} (and the fact that the winning parameter is independent of $\mathbf{b}$), the properties of Schmidt games (Section~\ref{subsecBWSCF}), Proposition~5.1 of~\cite{KW},\footnote{Thanks to Barak Weiss for pointing out this proposition.} and Theorem~3.1 of~\cite{F2} immediately imply the following corollary, which in particular gives  Kleinbock's main conjecture from~\cite{Kl}.  Note that an absolutely friendly measure is also absolutely decaying.  See Section~\ref{subsecBoFSoM} for details on these and on the fitting property of the measure $\mu$.

\begin{coro}\label{candidate-corollary}
Let $K \subset M_{m,n}(\RR)$ be the support of an absolutely friendly measure $\mu$.  Then, for any countable sequence $\{\mathbf{b}_i\} \subset  \RR^m$, \[K \cap (\cap_i \Ba^{\mathbf{b}_i}(m,n))\] is a winning set on $K$ and has Hausdorff dimension greater than or equal to $\underline{d}_\mu(K)$.  
If, in addition, $\mu$ is $\dim(K)$-fitting, then \[K \cap (\cap_i \Ba^{\mathbf{b}_i}(m,n))\] also has Hausdorff dimension equal to $\dim(K)$.
\end{coro}

\begin{rema}  In the corollary, if there exist constants $c_1, c_2, r_0 >0$ such that \[c_1 r^{\dim(K)} \leq \mu(B(x, r)) \leq c_2r^{\dim(K)},\] whenever $r\leq r_0$ and $x \in K$, then $\underline{d}_\mu(K) = \dim(K)$.
\end{rema}

Our second result is a generalization of the main result (Theorem 2) of~\cite{BHKV} to winning sets and to a larger class of fractals.  The result of~\cite{BHKV}, which shows full Hausdorff dimension, requires a high degree (related to $m$, see~\cite{BHKV} for the precise formulation) of regularity of the fractal.  This high degree of regularity precludes some common fractals (the Cantor set for example) that are included in Theorem~\ref{thmBAwinning} and \ref{thmDualBAwinning}.  In addition, Theorem~\ref{thmDualBAwinning} also generalizes the main result of~\cite{T2} to any dimension.  

\begin{theo} \label{thmDualBAwinning}
Let $K \subset \RR^m$ be the support of an absolutely $\eta$-decaying measure $\mu$.  Then, for any $A \in M_{m,n}(\RR)$, \[K \cap \Ba_A(m,n)\] is a winning set on $K$.
\end{theo}

 We again note that the winning parameter is a positive real number, independent of $A$. 

The proof of Theorem~\ref{thmDualBAwinning} in Section~\ref{secPTDualBAwinning}, which uses the space of unimodular lattices, is different from the second-named author's proof in~\cite{T2} of the special case $n=m=1$ and $K=\RR$, which uses continued fractions. For general $n$ and $m \in \NN$ and $K= \RR^m$, N.~Moshchevitin has a second proof that $\Ba_A(m,n)$ is winning for any $A$ which uses yet a third technique involving lacunary sequences~\cite{Mo}.  
 To our knowledge, Moshchevitin's remarkable proof, which is close to Schmidt's original proof that $\Ba^{\mathbf{0}}(m,n)$ is winning, does not give Theorem~\ref{thmDualBAwinning}.  Also, just before the finishing of the writing of this paper, we received the preprint~\cite{BFK} which gives an alternate proof of Theorem~\ref{thmDualBAwinning}.
 
  We would like to point out that U.~Shapira recently obtained a theorem concerning the set of multiplicative badly approximable systems, see \cite{Sha}.  In contrast to the results here there it is shown that for certain (and also almost all) $A\in {M}_{(1,2)}$ (resp.\ $A\in M_{(2,1)}$) the set of multiplicative badly approximable numbers $b\in\RR$ (resp.\ vectors $\mathbf{b}\in\RR^2$) can be empty.

Theorem~\ref{thmDualBAwinning}, the properties of Schmidt games (Section~\ref{subsecBWSCF}), and Theorem~3.1 of~\cite{F2} again immediately imply a corollary regarding intersections of  $\Ba_{A_i}(m,n))$, just as in Corollary~\ref{candidate-corollary}.

Finally, for singular matrices $A$, we can strengthen Theorem~\ref{thmDualBAwinning} by only considering the subset of infinitely badly approximable vectors $\Ba_A^\infty(,m,n)$ and obtain the following which is proven in Section~\ref{secAPthmDualBAwinningSing}.

\begin{theo} \label{thmDualBAwinningSing}
Let $K \subset \RR^m$ be the support of an absolutely $\eta$-decaying measure $\mu$.  Then, for any $A \in \mathbf{SM}_{m,n}(\RR)$, \[K \cap \Ba^\infty_A(m,n)\] is a winning set on $K$.
\end{theo}

A corollary like Corollary~\ref{candidate-corollary} also follows immediately.

We introduce winning sets and the space of unimodular lattices in Section~\ref{secBackground}, where we also introduce our method
in the classical case of $\mathbf{b}=\mathbf{0}$ and the Lebesuge measure. 
In Section~\ref{secPTDualBAwinning} we turn to a proof of our second result, Theorem~\ref{thmDualBAwinning}.  
In Section~\ref{secAPothmBAwinning} we prove Theorem~\ref{thmBAwinning} by showing how to extend the strategy in \cite{Sch3} resp.\ \cite{F1}.  Our third result, Theorem~\ref{thmDualBAwinningSing}, is a corollary of the proof of our second result and is presented in Section~\ref{secAPthmDualBAwinningSing}.

\section{Background}\label{secBackground}

The proofs of our results require two tools:  Schmidt games (see~\cite{Sch2} for a reference) and the basic concepts concerning flows on the space of unimodular lattices (see Chapter 9 of~\cite{EW} or~\cite{BM} for a reference).  In Section~\ref{subsecBWSCF}, we introduce the first tool, and, in Section~\ref{subsecBWSCF2}, we introduce the second.  
Finally our results are for fractals supported on certain measures, which we introduce in Section~\ref{subsecBoFSoM}.  (See, for example,~\cite{KLW}, \cite{F1}, and~\cite{PV} for additional details on these fractals.)  

\subsection{Schmidt games and winning sets}\label{subsecBWSCF}

W.~Schmidt introduced the games which now bear his name in~\cite{Sch2}.  Let $S$ be a subset of a complete metric space $M$.  For any point $x \in M$ and any $r \in \RR_+$, we denote the closed ball in $M$ around $x$ of radius $r$ by $B(x, r)$.  Even though it is possible for there to exist another $x' \in M$ and $r' \in \RR_+$ for which $B(x,r) = B(x',r')$ as sets in $M$, there will not be any ambiguity for us as we will always assume that we have chosen (either explicitly or implicitly) a center and a radius for each closed ball.  Let $\rho(A)$ denote the radius of the closed ball $A$. 
Schmidt games require two parameters:  $0 < \a <1$ and $0 < \b<1$.  Once values for the two parameters are chosen, we refer to the game as the $(\a,\b)$-game, which we now describe.  Two players, Black and White, alternate choosing nested closed balls $B_1 \supset W_1 \supset B_2 \supset W_2 \cdots$ on $M$ such that $\rho(W_n) = \a \rho(B_n)$ and $\rho(B_n) = \b \rho(W_{n-1})$.  The second player, White, \textit{wins} if the intersection of these balls lies in $S$.  A set $S$ is called \textit{$(\a, \b)$-winning} if White can always win for the given $\a$ and $\b$.  A set $S$ is called \textit{$\a$-winning} if White can always win for the given $\a$ and any $\b$; here $\a$ is called the \textit{winning parameter}.  A set $S$ is called \textit{winning} if it is $\a$-winning for some $\a$.  
Schmidt games have four important properties for us~\cite{Sch2}: \medskip


\noindent$\bullet$ Countable intersections of $\a$-winning sets are again $\a$-winning.




\noindent$\bullet$ Let $0 < \a \leq 1/2$.  If a set in a Banach space of positive dimension is $\a$-winning, then the set with a countable number of points removed is also $\a$-winning.

\noindent$\bullet$ The sets in $\RR^m$ which are $\a$-winning have full Hausdorff dimension.

Note that the last property has been generalized in two (related) ways.  Theorem~3.1 of~\cite{F2} states that, for a closed set $K \subset \RR^m$ which is the support of an absolutely $\eta$-friendly and $\dim(K)$-fitting measure, the $\a$-winning sets on $K$ have the same Hausdorff dimension as $K$.  Proposition~5.1 of~\cite{KW} states that, for $K$ the support of a Federer measure, the Hausdorff dimension of winning sets are greater than or equal to $\underline{d}_\mu(K)$.  See Section \ref{subsecBoFSoM} for definitions.

\subsection{The space of unimodular lattices}\label{subsecBWSCF2}
Let us now discuss a flow on the space of unimodular lattices and its relationship to systems of affine forms.  Let $\langle A, \mathbf{b} \rangle \in \widetilde{M}_{m,n}(\RR)$ and $k = m+n$.  The product $A \mathbf{q}$ can be viewed as a collection of $m$ linear forms in $n$ variables $q_1,\ldots,q_n$.  For non-zero $\mathbf{b}$, we call the expression $A \mathbf{q} - \mathbf{b}$ a \textit{system of $m$ affine forms in $n$ variables.}  We are interested in the size of $\|A \mathbf{q} - \mathbf{b}\|_{\ZZ}$ and $\|\mathbf{q}\|$ for $\mathbf{q} \in \ZZ^n$.  Let us combine all this data by considering the $(k+1) \times (k+1)$ matrix \[L_A(\mathbf{b}) := \begin{pmatrix}
   I_m & A & \mathbf{-b} \\
   0 & I_n & 0\\
   
   0 & 0 & 1
   \end{pmatrix},
  \] 
  where $I_\ell$ denotes the $\ell \times \ell$ identity matrix.  
 Moreover, we introduce the $k=m+n$-dimensional affine lattice 
\[
  L_A(\mathbf{b}) (\ZZ^{k}\times\{1\}):=\Bigg\{ L_A(\mathbf{b}) \begin{pmatrix}
   \mathbf{p} \\
   \mathbf{q} \\
   1
   \end{pmatrix}\mid \mathbf{p} \in \ZZ^m, \mathbf{q} \in \ZZ^n\Bigg\}
\] 
inside the ambient space $\RR^{k}\times\{1\}\cong\RR^{k}$. We will always identify $\RR^{k}$ with this affine subspace of $\RR^{k+1}$,
and will write $L_A(\mathbf{b}) (\ZZ^{k})$ as a shorthand for $L_A(\mathbf{b}) (\ZZ^{k}\times\{1\})$.
Finally, we define for any $t\in\RR$ the matrix
\[
 g_t := \begin{pmatrix}e^{t/m} I_m & 0 & 0\\ 0 & e^{-t/n} I_n& 0 \\ 0 & 0 & 1\end{pmatrix},
 \]
 which acts naturally on $\RR^{k+1}$ and also on $\RR^k$ (i.e.\ by the identification with the invariant affine subspace of $\RR^k\times\{1\}$). 
 The space $\Omega_{k,\mathrm{aff}}$ of \textit{affine unimodular lattices} in $k$-dimensions 
 is the space of all translates $\Lambda+\mathbf{c}\subset\RR^k$ 
 of unimodular lattices $\Lambda=g\ZZ^k$ for $g\in\operatorname{SL}(k,\RR)$ and $\mathbf{c}\in\RR^k$.
 All affine lattices $\Lambda+\mathbf{c}\subset\RR^k$ that we consider will be unimodular, and we often will think of them as subsets of $\RR^{k+1}$ in the way described above. In particular, the matrix $g_t$ acts on $\Omega_{k,\mathrm{aff}}$.


We call $\RR^k$ the \textit{time-particle space}.  When we refer to the \textit{origin} without further qualifications, we shall mean the \textit{origin of the time-particle space}.  We call $\{0\}^m\times\RR^n$ the \textit{time space} and $\RR^m\times\{0\}^n$ the \textit{particle space}.  The notions \textit{time component} and \textit{particle component} of a vector in $\RR^{k}$ are now clear. This terminology is explained by interpreting the elements $A\mathbf{q}+\ZZ^m\in \RR^m/\ZZ^m$ as the elements of the orbit of a $\ZZ^n$-action by rotation on the $m$-dimensional torus.
We let $\{\mathbf{e}_1, \cdots, \mathbf{e}_{n+m}\}$ denote the standard basis.

We will refer to $\Lambda$ as the associated lattice to the affine lattice $\Lambda+\mathbf{c}$. A subspace $V\subset\RR^k$ is called $\Lambda$-rational if $\Lambda\cap V$ spans $V$. 

For an $\ell$-dimensional parallelotope $P$, let $|P|$ denote its $\ell$-dimensional volume. If $V$ is a $\Lambda$-rational $\ell$-dimensional subspace we also write $|V|$ for the $\ell$-dimensional volume of the parallelotope $P\subset V$ spanned by a $\ZZ$-basis of $V\cap\Lambda$.  
A hyperplane $V$ (always of dimension $k-1$) is called \textit{small} if it is $\Lambda$-rational and $|V| \leq \xi_0:=\sqrt{k}$ and is called \textit{big} otherwise. 

All of the above notions are of course relative to an affine lattice $\Lambda+\textit{c}$. However, we will apply various elements of the flow $g_t$ to the affine lattice. In this case we will not always indicate this clearly, but if $H$ is $\Lambda$-rational and we talk about the covolume $|g_tH|$ then this is meant with respect to $g_t\Lambda$.  Furthermore, we say that a (big or small) hyperplane $H$ \textit{remains small} (with respect to $\Lambda+\textit{c}$) if there exists some $T_0\in\RR$ such that for all $t \geq T_0$, $ g_t H$ is small with respect to $g_t(\Lambda+\textit{c})$.

Also, we will use the following modification of a well-known theorem (Theorem 2.20 of~\cite{Da1}) due to S.~G.~Dani:

\begin{theo} \label{thmDani}
We have $\langle A,\mathbf{b}\rangle \in \Ba(m,n)$ if and only if all points in all affine lattices of the trajectory $\{g_t L_A(\mathbf{b}) \ZZ^k \mid t \in \RR_+\}$ are uniformly bounded away from the origin of the time-particle space.
\end{theo}

Even if the flow is replaced with a discrete time system by sampling times with uniformly bounded consecutive differences, the theorem still holds.
We also note that unlike the classical case of $\mathbf{b}=\mathbf{0}$, the above theorem does not relate the property $\langle A,\mathbf{b}\rangle \in \Ba(m,n)$ with the question whether the trajectory is bounded (i.e.\ has compact closure).

We now list a geometric lemma concerning the relationship between volume and unimodular lattices, which is straight forward to check.

\begin{lemm}\label{lemmDBPH}\label{lemmEVL1}
Let $\Lambda\subset\RR^k$ be a unimodular lattice. Let $H$ be a $\Lambda$-rational hyperplane. The distance between any two nearest parallel cosets $H+\mathbf{v}_1$ and $H+\mathbf{v}_2$ with $\mathbf{v}_1,\mathbf{v}_2\in\Lambda$ is equal to $1/ |H|$.
In particular, if the distance is $1/|H|<\xi_0^{-1}$ then the hyperplane $H$ is big.
 In any set of $k$ linearly independent vectors in $\Lambda$, there exists at least one lattice vector of length $\geq 1$.
\end{lemm}

Finally, we explain why small hyperplanes exist. The precise value of $\xi_0=\sqrt{k}$ is irrelevant for the main result of
the paper.  We also remark that for any unimodular lattice $\Lambda\subset\RR^k$ there exists only a finite number of small hyperplanes (but that this number cannot be bounded independent of the lattice). Both, the corollary regarding the existence of a small hyperplanes and the finiteness of the number of small hyperplanes follow from considering the dual lattice. Here the dual of a lattice $\Lambda\subset\RR^k$ is defined by $\Lambda^*=\{w\in\RR^k: \langle v,w\rangle\in\ZZ$ for all $v\in\Lambda\}$ and we note that there is a correspondence between primitive vector $w\in\Lambda^*$ and hyperplane $H=w^\perp$ for which $|H|$ with respect to $\Lambda$ equals $\|w\|$.

\begin{rema}\label{rmkAVHGZoI}
We also remark that the asymptotic volume of any hyperplane goes to  either zero or infinity in the following sense. Let $H$ be a $\Lambda$-rational hyperplane, then either $|g_tH|$ measured with respect to $g_t\Lambda$ goes to infinity or to $0$ as $t\to\infty$.\footnote{This is precisely the behavior that is also explained by considering the eigenvalues of $\bigwedge^{k-1}g_t$ acting on $\bigwedge^{k-1}\RR^k$, which leads to a formal proof.} To see this assume first that  $H$ contains the time space $\{0\}^m\times\RR^n$. In this case $H$ is spanned by the time space and a hyperplane of the particle space, is invariant under $g_t$ and $g_t$ restricted to $H$ has determinant $e^{-t/m}$. This shows clearly that $|g_tH|$ with respect to $g_t\Lambda$ goes to zero. In the second case $H$ is spanned by $m$ vectors that project to a basis of the particle space $\RR^m\times\{0\}^n$ and by $n-1$ vectors that belong to the time space. In this case it follows that $|g_tH|$ measured with respect to $g_t\Lambda$ goes to infinity.
\end{rema}


\subsection{Fractals supported on measures}\label{subsecBoFSoM}

Let $\mathcal{L}$ denote
an affine $(n-1)$-dimensional hyperplane of $\RR^n$. For $\epsilon > 0$, let
$\mathcal{L}^{(\epsilon )}$ denote the $\epsilon$-thickening of
$\mathcal{L}$.  A locally finite Borel measure $\mu$ on $\RR^n$ is called \textit{absolutely
$\eta$-decaying} if there exist
   strictly positive constants $C, \eta$ and $r_0$ such that for any
    hyperplane $\mathcal{L}$, any $\epsilon >0$, any $x \in \textrm{supp}(\mu)$,
    and any positive $r< r_0$,
    \begin{equation*}
      \mu(B(x, r) \cap \mathcal{L}^{(\epsilon )}) \leq C \left(
        \frac{\epsilon}{r} \right)^{\eta}
      \mu(B(x, r)).
    \end{equation*}  A locally finite Borel measure $\mu$ is called \textit{Federer (or doubling)} if there exist strictly positive constants $D$ and $r_0$ such that, for any $x \in \textrm{supp}(\mu)$ and any positive $r< r_0$, 
   \[
     \mu(B(x,\frac{1}{2} r))> D\mu(B(x ,r)).
   \]  An absolutely $\eta$-decaying, Federer measure $\mu$ is called \textit{absolutely $\eta$-friendly}.  
    
For a metric space $(X, d)$, a given $x \in X$ , and real numbers $r > 0, 0 < \b < 1,$ let $N_X (\b , x, r)$ denote (following~\cite{F2}) the maximum number of disjoint balls (centered at a point of $X$) of radius $\b r$ contained in $B(x, r)$.  A locally finite Borel measure $\mu$ is \textit{$\delta$-fitting} if there exist constants $0<r_1 \leq 1, M , \textrm{ and } \delta$ such that, for every $0 < r \leq r_1 , 0 < \b < 1$ and $x \in \textrm{supp}(\mu)$, \[
 N_{\textrm{supp}(\mu)} (\b , x, r) \geq M \b^{-\delta}.
\]  
    
Lebesgue measure on $\RR^n$ is an example of an absolutely friendly, fitting measure.  Besides $\RR^n$, the support of an absolutely friendly, fitting measure includes the Cantor set, the Koch curve, the Sierpinski gasket, or, in general, the attractor of 
an irreducible finite family of contracting similarity maps of $\RR^n$ satisfying the open set condition (see Corollary~5.3 of~\cite{F2} and Theorem~2.3 of~\cite{KLW} for more details).

\section{Proof of Theorem~\ref{thmDualBAwinning}} \label{secPTDualBAwinning}

In this section, we prove Theorem~\ref{thmDualBAwinning}.  An understanding of this proof will illuminate the proofs of our other results.  
The proof consists in describing the strategy that  player White should use, and in proving that White indeed always wins by using this strategy.
Note that the matrix $A$ and so the lattice $\Lambda=L_A(0)\ZZ^k$ are given by assumption while the game takes place on the set of possible translations $\mathbf{b}$ which define the affine lattices $\Lambda-\begin{pmatrix}\mathbf{b}\\ 0\end{pmatrix}=L_A(\mathbf{b})\ZZ^k$. 

Let $0<\b<1$ be fixed, and note that $\eta$ and $C$ are two constants coming from the definition of absolute $\eta$-decay, which we assume for $\mu$. By our assumption, $K=\operatorname{supp}\mu$.
Let 
\begin{align*}
\a &< \bigl(4(2\xi_0C)^{\frac1\eta}\bigr)^{-1}\\
T &= -m\log(\a\b).
\end{align*}  
Our strategy will use the value of $\b$ implicitly by using the transformation $g_T$ on $\RR^k$. Also note that $\alpha$ has been chosen independent of $\b$ (which is required for showing that the game is $\alpha$-winning).

Let us point out the crucial link between steps of the game and applications of $g_T$. In every complete cycle of the game, the radii of the  balls $B_\ell\supset W_\ell$ are multiplied by $\alpha\beta$ and the game then continues with the shrinked balls. In the dynamical system, we instead replace the given affine lattice $\Lambda_\ell$ (representing a point in $\Omega_{k,\mathrm{aff}}$) by the lattice $g_T\Lambda_\ell=\Lambda_{\ell+1}$. By definition the map $g_T$ expands the particle space by $(\alpha\beta)^{-1}$ and the time space is contracted (by $(\alpha\beta)^{\frac{m}n}$). Roughly speaking, this allows one to relate statements about the lattice $g_T^\ell L_A(\mathbf{b})\ZZ^k$ with respect to the unit ball to statements about elements of the $(\alpha\beta)^\ell$-ball in particle space and elements $\mathbf{q}\in\ZZ^n$ of the time space of size less than $(\alpha\beta)^{-\frac{m}n\ell}$ --- this is the basis of Theorem~\ref{thmDani}. White tries to restrict the choice of $\mathbf{b}$ by choosing the new ball (in the game of radius $\alpha\rho(B_\ell)$ and in the dynamical picture of radius $\alpha$) so that $g_T^{\ell+1}L_A(\mathbf{b})\ZZ^k$ has no elements in a ball around zero of some fixed radius independent of how $\mathbf{b}$ is chosen from the new ball. There is one potential problem in this simple-minded strategy, namely it could happen that the affine lattice $g_T^{\ell}L_A(\mathbf{b})\ZZ^k$ contains an $m$-dimensional subspace  that is close to the particle space $\RR^m\times\{0\}^n$ and on which the lattice points of $g_T^{\ell}L_A(0)\ZZ^k$ are highly dense and the center $\mathbf{b}$ at that stage is such that the affine lattice $g_T^{\ell}L_A(\mathbf{b}))\ZZ^k$ contains lattice elements in the unit ball. In this case, the lattice $g_T^{\ell+1}L_A(\mathbf{b})\ZZ^k$ will contain points close to zero independently of how $\mathbf{b}$ is chosen from $B_\ell$.  If $A$ is badly approximable itself, then this problem does not appear (as the lattices $g_T^\ell L_A(0)\ZZ^k$ for $\ell=1,2,\ldots$ remain uniformly discrete) and the strategy is quite straightforward. In general, the strategy of White is to study the behavior of rational hyperplanes and, by making correct moves earlier on in the game, the above bad scenario can be avoided by moving away from a hyperplane before it becomes very short. The assumption that $K$ supports an absolutely decaying measure is precisely the condition that allows White to move away from hyperplanes.

Also useful will be the following identities which formalizes some of the above discussions. First the affine lattice $L_A(\mathbf{b})\ZZ^k$ can be obtained from $L_A(0)\ZZ^k$ by application of the translation operator $L_0(\mathbf{b})$ since $L_0(\mathbf{b})L_A(0)=L_A(\mathbf{b})$. Second, application of $g_t$ to $L_A(\mathbf{b})\ZZ^k$ gives the same as application of the translation operation $L_0(e^{t/m}\mathbf{b})$ to $g_tL_A(0)\ZZ^k$ as
\[
 g_tL_0(\mathbf{b})g_t^{-1}=L_0(e^{t/m}\mathbf{b}).
\]

We continue with a formal description of the strategy. 
Depending on $A$ there are two cases; we begin with an easy but atypical case\footnote{This case is actually trivial as the $\ZZ^n$-orbit defined by $A$ on $\mathbb{T}^m$ is not even dense, but we give the simplified version of the argument used in the general case to show concretely why hyperplanes can be helpful.}.

\subsection{Case 1: There is an $L_A(0)\ZZ^k$-rational hyperplane whose covolume goes to zero.}
Suppose the ball $B_1=B(\mathbf{b}_1,\rho_1)\subset\RR^m$ with center $\mathbf{b}_1\in K$ and radius $\rho_1>0$ has been chosen by player Black. Let $H\subset\RR^k$ be the hyperplane for which $|g_tH|$ measured with respect to $g_tL_A(0)\ZZ^k$ goes to zero as $t\to\infty$.
As discussed in Remark \ref{rmkAVHGZoI}
this means that $H$ contains $\{0\}^m\times\RR^n$ and intersects the particle space $\RR^m\times\{0\}^n$ in a
 hyperplane. 
 We choose $t_0>0$ such that  $\rho_1=e^{-t_0/m}$. Moreover, we may assume that $g_{t_0}H$ is short, 
 in fact with covolume less than $\frac13$, with respect
 to the lattice $g_{t_0}L_A(0)\ZZ^k$. Otherwise we let White play a few steps without any particular goal other than making the balls smaller and  the corresponding parameter $t_0$ larger. 
Assuming now that the covolume of $g_{t_0}H$ w.r.t.\ $g_{t_0}L_A(0)\ZZ^k$ is less than $\frac13$, we see that distinct cosets $\mathbf{v}+g_{t_0}H$ for $\mathbf{v}\in g_{t_0}L_A(0)\ZZ^k$ need to be at least $3$ far apart. White wants to make sure that the element $\mathbf{b}$ constructed by the game is such that $L_A(\mathbf{b})\ZZ^k+H$ does not contain the origin. (In the case considered below, we will have to be more careful about the distance to such hyperplanes.) Assume that the coset $\mathbf{v}+g_{t_0}H$ for some $\mathbf{v}\in g_{t_0}L_A(0)\ZZ^k$ indeed intersects $e^{t_0/m}B(\mathbf{b}_1,\rho_1)$ --- by the distance of these cosets from one another
 there can be only one. Let $\mathcal{L}\subset \RR^m\times\{0\}^n$  be the hyperplane such that $g_{t_0}\mathcal{L}$ is the intersection of the coset $\mathbf{v}+g_{t_0}H$ with $\RR^m\times\{0\}^n$.
Applying the definition of absolutely decaying to the $\epsilon$-neighborhood $\mathcal{L}^{(\epsilon)}$ with $\epsilon=2\alpha\rho_1$ and the ball $B(\mathbf{b}_1,\rho_1(1-\alpha))$ it follows from the choice of $\alpha$ that there is some $\mathbf{b}_1'\in K\cap
B(\mathbf{b}_1,\rho_1(1-\alpha))\setminus \mathcal{L}^{(\epsilon)}$. The strategy of White is to choose one such point as the center of $W_1$ (which is allowed as $W_1\subset B_1=B(\mathbf{b}_1,\rho_1)$. After this first step
White does not have to be careful --- we claim that White wins independently of the remaining steps of the game. The reason for this is simply that the constructed $\mathbf{b}$ from the game must have $\mathbf{b}\notin\mathcal{L}^{(\alpha\rho_1)}$. This implies that $\begin{pmatrix}\mathbf{b}\\0\end{pmatrix}$ together with a basis of  $H\cap L_A(0)\ZZ^k$ span a parallelepiped of positive $k$-dimensional volume. As $g_t$ does not change the volume and the volume of the base of the parallelepiped inside $H$ goes to zero (as it equals the covolume of $H$), it follows that the distance of $g_t\begin{pmatrix}\mathbf{b}\\0\end{pmatrix}$ to $g_tH$ goes to infinity. The same applies to any other cosets of $H$, which shows that $g_tL_A(\mathbf{b})\ZZ^k$ can indeed not contain small vectors as $t\to\infty$. 
This concludes the proof of this simple case by Theorem \ref{thmDani}.

\subsection{Case 2:  No hyperplane of $L_A(0)\ZZ^k$ remains small.}

Let $B_1=B(\mathbf{b}_1,\rho_1)$ be the ball chosen by player Black. We define $t_1$ such that $e^{t_1/m}\rho_1=1$ and also the affine lattice $x_1=g_{t_1}L_A(\mathbf{b}_1)\ZZ^k$. We use induction to describe the strategy and the proof.  In the initial step of the induction, we ignore any (probably ridiculously) small hyperplanes of $g_{t_1} L_A(\mathbf{0})\ZZ^k$ and let White play without any strategy.  In later steps of the induction White will make sure that any small hyperplances $g_{t_1+(J-1)T}H$ have their cosets $v+ g_{t_1+(J-1)T}H$ for $v\in x_J=g_{t_1+(J-1)T}(L_A(\mathbf{b})\ZZ^k)$ at a significant distance from the origin. To simplify notation we define $t_J=t_1+(J-1)T$.

Since a small hyperplane always exists and since, in this case, no hyperplane remains small forever, at some future point, a big hyperplane must become small.  Let $J\geq 1$ be minimal such that there is a hyperplane $H$ such that $g_{t_{J-1}}H$ is big (w.r.t.\ $g_{t_{J-1}T}L_A(0)\ZZ^k$) but $g_{t_J}H$ is small (w.r.t.\ $g_{t_J}L_A(0)\ZZ^k$).  If there is more than one such hyperplanes, we choose $H$ such that $g_{t_k}H$ is small the longest (i.e.\ for the most $k>J$). White may play without any particular goal up to stage $J$ of the game. Suppose Black has chosen his ball $B_J=B(\mathbf{b}_J,\rho_J)$.  Consequently, we note that $\rho(e^{t_J/m} B_J) = 1$. 
This means that White is given the lattice $x_J=g_{t_J}L_A(\mathbf{b}_J)\ZZ^k$ and the freedom to replace $x_J$ by $L_0(\mathbf{b})x_J$ for any $\mathbf{b}\in B(0,1-\alpha)$. More precisely, this corresponds to choosing the center $\mathbf{b}_J+e^{-t_J/m}\mathbf{b}$ for the ball $W_J$ and White also has to ensure that this center belongs to $K$. 

Note that the hyperplane $H$ cannot contain the particle space $\RR^m\times\{0\}^n$ as otherwise the covolume of $H$ would be monotonically increasing (contradicting our reasons to look at $H$ in the first place). Moreover, we claim that the angle between $H$ and the particle space  $\RR^m\times\{0\}^n$ is significant in the following sense: There exists some $\delta>0$ (which depends on $k$ and $T$) such that for any vector $v\in\RR^m\times\{0\}^n$ which is in distance $d$ from $H\cap (\RR^m\times\{0\}^n)$ produces together with the $k-1$-dimensional parallelepiped in $g_{T_J}H$ corresponding to $y_J$, a $k$-dimensional parallelepiped of volume $\delta d |g_{T_J}H|$, where $|g_{T_J}H|$ denotes the $k-1$-dimensional volume of the parallelepiped in $g_{T_J}H$.

To see the existence of $\delta$, recall that $|g_{T_J}H|$ equals the norm of the vector $v_1\wedge \cdots \wedge v_{k-1}$ where $v_1,\ldots, v_{k-1}$ is a basis of $g_{T_J}H\cap y_J$. Furthermore, $\bigwedge^{k-1} g_T$ has eigenvalues $e^{-T/m}$ of multiplicity $m$ (corresponding to those hyperplanes that contain the time space) and $e^{T/n}$ of multiplicity $n$ (corresponding to those hyperplanes that contain the particle space). As $g_JH$ is big but $g_{J+1}H$ is small the vector $v_1\wedge\ldots\wedge v_{k-1}$ splits into a sum of eigenvectors $w_-$ with eigenvalue $e^{-T/m}$ and $w_+$ with eigenvalue $e^{T/n}$. A simple calculus exercise now shows that since the size of the vector decreases from $J-1$ to $J$ the vector $w_-$ must be significant and cannot be much smaller than $w_+$. 
Finally when calculating the volume $|v_1\wedge\cdots\wedge v_{k-1}\wedge v|$ of the $k$-dimensional parallelepiped mentioned above, the component $w^+$ is irrelevant as $w_+\wedge v=0$. This gives the claim.

The covolume of $g_{t_J}H$ is $\leq \xi_0$, and so the distance between any two cosets of elements in $x_J$ with respect to $g_{t_J}(H)$ must be $\geq \xi_0^{-1}$. This implies that at most $2\xi_0$-many of the cosets $v+g_{t_J}H$ with $v\in x_J$ which can intersect the unit ball. Taking those intersections into account, the strategy of White is such to put his new ball $B=B(\mathbf{c},\a)\subset B(0,1)$, in the dynamical picture, with center $\mathbf{c}\in B(0,1-\alpha)$ such that after the shift $L_0(\mathbf{b})$ by any $\mathbf{b}\in B$
the distance of $(L_0(\mathbf{b})x_J+g_{t_J}H)\cap(\RR^m\times\{0\}^n)$ to the origin is at least $\alpha$. Note that this intersection
consists of cosets of $\mathcal{L}$ of which there are at most $2\xi_0$ many which are in danger of getting, after the shift, close to the origin.

Of course, White is obliged to 
make his choice of $\mathbf{c}$ such that $\mathbf{b}_J+e^{-t_J/m}\mathbf{c}$, namely the center of the ball of White in the game, also belongs to $K$. We have choosen $\alpha$ in such a way that after applying the condition of absolute $\eta$-decay $2\xi_0$-many times for $\epsilon=2\alpha\rho_J$ and $r=\rho_J(1-\alpha)\geq \frac12\rho_J$ we are still ensured to find an element of $K$ outside the $2\alpha\rho_J$-neighborhoods of the $\leq 2\xi_0$ cosets of the hyperplane $\mathcal{L}$ that are relevant. 

The above strategy ensures that the volume of  the $k$-dimensional pyramid that is spanned by the $k-1$-dimensional parallelepiped (with $k-1$-dimensional volume $|g_{t_J}H|\geq \xi_0e^{-T/m}$) inside any of the cosets of $v+g_{t_J}H$ with $v\in L_0(\mathbf{b})x_J$ and $\mathbf{b}\in B$ has volume at least $\xi_0e^{-T/m}\alpha\delta$ (which we agree to call significant as it doesn't depend on $J$). Since $g_T$ does not change this volume, 
we see that the smallest vector of $g_{kT}(L_0(\mathbf{b})x_J)$ has norm at least 
\begin{equation}\label{smallest-norm}
\xi_0e^{-T/m}\alpha\delta  |g_{t_k}H|^{-1}\geq  e^{-T/m}\alpha\delta,
\end{equation}
where the last inequality holds for any $k\geq 0$ with $|g_{t_{J+k}}H|\leq \xi_0$. For those times, White is protected from getting short vectors in the corresponding affine lattices.

 If for some $J'>J$
there is another hyperplane $H'$ that just became small as time $J'$, then White has to repeat the above procedure, again playing to make the volumes of certain pyramids significant. This may and eventually will add protection time. Repeating the procedure infinitely often constructs some shift $\mathbf{b}_\infty$. The construction (the protection times cover in the end the interval $[J,\infty)$) and  Theorem \ref{thmDani} imply $(A,\mathbf{b}_\infty)\in\Ba(m,n)$. 

\subsection{Proof of Theorem~\ref{thmDualBAwinningSing}}\label{secAPthmDualBAwinningSing} 

As mentioned in the introduction, Theorem~\ref{thmDualBAwinningSing} is really a corollary of the proof of Theorem~\ref{thmDualBAwinning} in the following sense. Assume now that $A$ is singular, and let White use the same strategy as described above. Then after the game has finished, it has constructed some $\mathbf{b}_\infty$. Let $x=L_A(\mathbf{b}_\infty)\ZZ^k$
be the corresponding affine lattice. Then for large enough $J$ we will have by \eqref{smallest-norm} that $g_k (x)$ has no vector that is shorter than $\xi_0e^{-T/m}\alpha\delta  |g_{t_k}H|^{-1}$ where $H$ is the hyperplane for which $|g_{t_k}H|$ is smallest.
However, if $A$ is singular, then applying the Mahler compactness criterion to the dual lattice, it follows that 
$ \min_H |g_{t_k}H|$ goes to zero. Therefore, the norm of the smallest element of $g_t(x)$ goes to infinity. This implies that $\mathbf{b}_\infty\in\Ba_A^\infty(m,n)$ by (a simple strengthening of) Theorem \ref{thmDani}.

\section{Proof of Theorem~\ref{thmBAwinning}} \label{secAPothmBAwinning} 

We now show how to adapt the available strategies for White and $\Ba^{0}(m,n)$ in \cite{Sch3} for the Lebesgue measures, resp.\ in \cite{F1} for friendly measures, to get a strategy for $\Ba^{\mathbf{b}}(m,n)$.  So let $\alpha_0$ be a value so that 
$\Ba^{0}(m,n)$ is $(\alpha_0,\beta_0)$-winning for every $\beta_0>0$. We define $\alpha=\bigl(4(2C)^{\frac1\eta}\bigr)^{-1}\alpha_0$, and for every $\beta>0$ we define $\beta_0=\beta\bigl(4(2C)^{\frac1\eta}\bigr)^{-1}$ so that $\alpha\beta=\alpha_0\beta_0$. The strategy is, for every given ball $B_\ell=B(A_\ell,\rho_\ell)$ to use the known strategy of White for $\Ba^{0}(m,n)$
to choose $B_\ell'=B(A_\ell',\alpha_0\rho_\ell)\subset B_\ell$ with $A_\ell'\in K$ and an additional step below for $\Ba^{\mathbf{b}}(m,n)$
to get a ball $W_\ell\subset B_\ell'$ of radius $\alpha\rho_\ell$ and center in $K$. We then show that this modified strategy is winning for   $\Ba^{\mathbf{b}}(m,n)$. 

We may and will assume $\mathbf{b}\notin\ZZ^m$. Given $B_\ell'=B(A_\ell',\alpha_0\rho_\ell)$ we define the affine lattice $x_\ell=g_{t_\ell}L_{A_\ell'}(\mathbf{b})\ZZ^k$. Here $t_\ell$ is chosen such that $g_{t_\ell}L_D(0)=L_{  (\alpha_0\rho_\ell)^{-1} D }(0)g_{t_\ell}$ for any $D \in M_{n,m}(\RR)$. As in the argument above, this makes the additional step of choosing the subball $W_\ell$ and replacing $x_\ell$ with the lattice corresponding to the new center equivalent to choosing a subball of $B(0,1)$ of radius $\bigl(4(2C)^{\frac1\eta}\bigr)^{-1}$ and applying the center to $x_\ell$.

Let $\mathbf{v}\in x_\ell$ be a vector of smallest norm. 
We are choosing the new center in such a way that the particle component $\mathbf{v}_p$ is significant in relationship to the norm $\|\mathbf{v}_t\|$ of the time component.\footnote{Henceforth, we use the subscript $t$ on time-particle vectors to denote their time components and use $p$ to denote particle components.}
Indeed, there is a proper affine subspace $\mathcal L\subset M_{m,n}$ (which depends on $\mathbf{v}$)
such that $(L_A(0)\mathbf{v})_p=0$ if and only if $A\in\mathcal L$.\footnote{Here $\mathbf{v}_t \neq \mathbf{0}$.  If $\mathbf{v}_t = \mathbf{0}$, then the calculation in the rest of this paragraph is trivial since $L_D(0)$ fixes $\mathbf{v} = \mathbf{v}_p$ for any $D$.} Moreover, it is straight forward to check that $A\notin \mathcal L^{(\varepsilon)}$ implies $\|(L_A(0)\mathbf{v})_p\|\geq\varepsilon\|\mathbf{v}_t\|$.
By definition of absolute $\eta$-decay, applied to $B(A_\ell',\frac12\alpha_0\rho_\ell)$ and $\varepsilon=2\bigl(4(2C)^{\frac1\eta}\bigr)^{-1}\alpha_0\rho_\ell$, we are sure to find a new center $A_\ell''\in K\setminus\mathcal L^{(\varepsilon)} $
with $B(A_\ell'',\alpha\rho_\ell)\subset B(A_\ell',\alpha_0\rho_\ell)$ corresponding (in the sense described above and depending on $x_\ell$) to a subball $B$ of $B(0,1)$ such that 
\begin{equation}\label{short-vector-change}
   \|(L_D(0)\mathbf{v})_p\|\geq \bigl(4(2C)^{\frac1\eta}\bigr)^{-1} \|\mathbf{v}_t\|
\end{equation}
whenever $D\in B$. 

We now prove that the above strategy for White is winning. By the assumed strategy, the matrix $A$ that belongs to the intersection of all balls is badly approximable. So let $\epsilon>0$ be small enough so that $g_tL_A(0)\ZZ^k$ does not contain any nonzero element of norm $\leq\epsilon$ for any $t\geq 0$. We may also choose $\delta>0$ such that the affine lattice $L_A(\mathbf{b})\ZZ^k$ does not contain any element of $B(0,\delta)$ (as $\mathbf{b}\notin\ZZ^k$). Finally suppose $c>0$ is such that $\|D\mathbf{v}_t\|\leq c\|\mathbf{v}_t\|$ for all $D\in B(0,1)$ and $\mathbf{v}\in\RR^k$. Then we claim that $g_{t_\ell}L_A(\mathbf{b})\ZZ^k$ does not contain any element of $B\bigl(0,r)$
for
$$
r=\min (\delta,\epsilon) (1+c)^{-2}(\alpha\beta)^{\frac{n}{m+n}}/2
$$
and for any $\ell\geq 0$. By Theorem~\ref{thmDani} this claim implies that $A\in\Ba^{\mathbf{b}}(m,n)$.    

The claim holds for $\ell=0$ by choice of $\delta$. Now suppose $\mathbf{w}_0\in\bigl(g_{t_\ell}L_A(\mathbf{b})\ZZ^k\bigr)\cap B(0,r)$ exists and $\ell\geq 1$ is chosen minimally with this property. Then the affine lattice $x_\ell=g_{t_\ell}L_{A_\ell'}(\mathbf{b})\ZZ^k$ that was used in the strategy differs from $g_{t_\ell}L_{A}(\mathbf{b})\ZZ^k$ by an application of $L_D(0)$ with some $D\in B(0,1)$ --- just because $A$ belongs to the ball that was chosen by White at stage $\ell$. Therefore, there exists a vector $\mathbf{w}\in x_\ell\cap B(0,(1+c)r)$. Going back one step in the dynamical iteration corresponding to the game
we get $x_{\ell}=g_TL_{D'}(0)x_{\ell-1}$ where $D'\in B(0,1)$. Here $T$ is such that $L_A(0)g_T=g_TL_{(\alpha\beta)^{-1} A}(0)$, i.e.\ $g_T$ has eigenvalues $(\alpha\beta)^{-\frac{n}{m+n}}$ and $(\alpha\beta)^{\frac{m}{m+n}}$. Therefore, $x_{\ell-1}$ contains an element $\mathbf{v}'$ of norm $\leq (1+c)^2(\alpha\beta)^{-\frac{n}{m+n}}r\leq \epsilon/2$. However, as $g_{t_{\ell-1}}L_A(0)\ZZ^k$ does not contain any nonzero element of norm $\leq \epsilon$ this shows that $\mathbf{v}'=\mathbf{v}$ is the element that was used in the additional step of the strategy at step $\ell-1$.\footnote{Since $g_{t_{\ell-1}}L_A(\mathbf{b})\ZZ^k$ is just a translation along a direction in particle space of $g_{t_{\ell-1}}L_A(0)\ZZ^k$, an $\epsilon/2$-ball can contain at most one lattice point; thus $\mathbf{v}' = \mathbf{v}$, the smallest vector of $x_{\ell-1}$.} Therefore, \eqref{short-vector-change} holds for $\mathbf{v}$.
Moreover, as $g_{t_{\ell-1}}\bigl(L_A(\mathbf{b})\ZZ^k\bigr)$ does not contain an element of norm less than $r$
we see that $\|\mathbf{v}\|\geq (1+c)^{-1}r$.\footnote{Note that $L_{D'}(0)$ fixes $\mathbf{v}_t$; moreover, its effect on $\mathbf{v}_p$ is small if $\mathbf{v}_t$ is small.} Therefore, one can derive (e.g.\ by considering the case
 $c\|\mathbf{v}_t\|\leq \frac12\|\mathbf{v}_p\|$ and the case  $c\|\mathbf{v}_t\|\geq \frac12\|\mathbf{v}_p\|$ separately) from \eqref{short-vector-change} that $(L_{D'}(0)\mathbf{v})_p$ is of size $\geq \kappa r$, where the constant $\kappa$ depends on $C,\eta,k,c$. As the particle space gets uniformly epxanded, this implies that after applying $g_T$ we have that $\mathbf{w}$ has norm $\|\mathbf{w}\|\geq\kappa r (\alpha\beta)^{-\frac{n}{n+m}}$. On the other hand we already know that $\|\mathbf{w}\|\leq (1+c)r$, which gives $\kappa r (\alpha\beta)^{-\frac{n}{n+m}}\leq (1+c)r$. This is
  a contradiction to the assumption that the claim does not hold,  if only $\beta$ is sufficiently small. Note that for Schmidt games it is allowed to assume that $\beta$ is sufficiently small  --- if White decides to use his strategy only every $p$ step of the game, this has the effect of replacing $\beta$ by the much smaller $\beta(\alpha\beta)^{p-1}$.

\section{Conclusion}


A badly approximable system of affine forms is determined by a matrix and a vector.  Our two main results, Theorems~\ref{thmBAwinning} and~\ref{thmDualBAwinning}, determine the size of the set of badly approximable systems of affine forms for a fixed vector and a fixed matrix respectively, and Theorem~\ref{thmBAwinning}, in particular, shows a fundamental conjecture in singly metric inhomogeneous number theory on the Hausdorff dimension of these sets for fixed vectors.  Moreover, our two theorems lead to another conjecture.  Instead of fixing either the vector or the matrix, one fixes neither and considers the size of $\Ba(m,n)$, which, recall, is the full set of badly approximable systems of affine forms.  A classical result, the doubly metric inhomogeneous Khintchine-Groshev Theorem (see Theorem II in Chapter VII of \cite{Ca} for the statement of the theorem), immediately implies that this set has zero Lebesgue measure.  With regard to dimension, Kleinbock has shown, using mixing of flows on the space of unimodular lattices, that the set has full Hausdorff dimension~\cite{Kl}.\footnote{It is interesting to note that, using the Marstrand slicing theorem (\cite{Fal}, Theorem 5.8), the main result in~\cite{BHKV}, Theorem~\ref{thmDualBAwinning}, the main result in~\cite{Mo}, Theorem~\ref{thmBAwinning}, and the main result in~\cite{BFK} give five additional, different proofs of Kleinbock's result.  Also, this result of Kleinbock immediately implies his other result mentioned in Section~\ref{secIntro} (see~\cite{Kl}).}  Moreover, Kleinbock conjectured that the set is winning~\cite{Kl} (or winning in the modified sense of~\cite{KW}, as mentioned in personal communication).  It seems interesting, although the authors have not yet undertaken this endeavor, to combine our proofs of Theorems~\ref{thmBAwinning} and~\ref{thmDualBAwinning} to yield a proof of not only this conjecture, but also a more general conjecture:  if $K \subset \widetilde{M}_{m,n}(\RR)$ is a closed subset supporting an absolutely $\eta$-decaying measure $\mu$, then $K \cap \Ba(m,n)$ is a winning set on $K$.  However, the obstacle to this could be the different ways in which $\mathbf{b}$ resp.\ $A$ in $L_A(\mathbf{b})$ are affected by conjugation with $g_t$; thus modified winning may be the better conjecture. 

\subsection{Strong winning}  Finally, we remark that the notion of strong winning for subsets of $\RR^n$ has recently been defined in \cite{Mc}.  Strong winning implies winning and is preserved by quasisymmetric homeomorphisms~\cite{Mc}.  It is not difficult to see that we can also conclude strong winning in Theorems~\ref{thmBAwinning},~\ref{thmDualBAwinning},~and~\ref{thmDualBAwinningSing} above.

\end{document}